\documentclass[12pt]{article}
\usepackage{amssymb,latexsym,amsmath}
\def\C{\mathbb C}
\def\R{\mathbb R}
\def\N{\mathbb N}
\def\Z{\mathbb Z}

\numberwithin{equation}{section}

\newtheorem{thm}{Theorem}[section]
\newtheorem{lem}{Lemma}[section]

\newtheorem{prop}{Proposition}[section]
\newtheorem{cor}{Corollary}[section]
\begin{document}
\sffamily

\title{Transcendental singularities for a meromorphic function with logarithmic derivative of finite lower order }
\author{J.K. Langley}
\maketitle

\begin{abstract}
It is shown that two key
results on transcendental singularities for meromorphic functions of finite lower order  have refinements
which hold under the weaker hypothesis
that  the logarithmic derivative
has finite lower order. 
\end{abstract}

\section{Introduction and results}

Suppose that $f$ is a  transcendental meromorphic function on $\C$ such that, 
as $z$ tends to infinity
along a path $\gamma $ in the plane, $f(z)$ tends to some  $\alpha \in \C$.
Then,  for each  $t > 0$, an unbounded subpath of 
$\gamma$ lies in a component 
$C(t)$ of the set $ \{ z \in \C : |f(z) -  \alpha | < t \}$.
Here $C(t) \subseteq C(s)$ if $0 < t < s$, and the intersection $\bigcap_{t>0} C(t)$ is empty \cite{BE}.  
The path $\gamma$ then determines a transcendental singularity of the inverse function $f^{-1}$  over the asymptotic value $\alpha$ and
each $C(t)$ is called a neighbourhood of the singularity  \cite{BE,Nev}. Two transcendental singularities over $\alpha$ are distinct if they have disjoint
neighbourhoods for some $t > 0$. 
Following \cite{BE,Nev}, a transcendental singularity of $f^{-1}$ over $\alpha  \in  \C$ is said to be direct
if $C(t)$, for some $t > 0$, contains finitely many points $z$ with  $f(z) = \alpha$, in which case there exists $t_1 > 0$
such that $C(t)$ contains no $\alpha$-points of $f$ for $0 < t < t_1$. 
A direct singularity over $\alpha \in \C$ is logarithmic if there exists $t > 0$ such that $\log t/(f(z)-\alpha)$ maps $C(t)$ 
conformally onto the right half plane. If, on the other hand, 
$C(t)$ contains infinitely many $\alpha$-points  of $f$, for every $t > 0$, then the singularity is called indirect:
a well known  example  is given by $f(z) = z^{-1}\sin z$, with $\alpha=0$ and $\gamma$  the positive
real axis $\R^+$.
Transcendental singularities of $f^{-1}$ over $\infty$ and their corresponding neighbourhoods 
may be defined and classified using $1/f$, and the asymptotic and critical values of $f$ together comprise the singular values of $f^{-1}$. 

If $f$ has finite (lower) order of growth, as defined in terms of the Nevanlinna characteristic function $T(r, f)$ \cite{Hay2,Nev}, 
then the number of direct singularities is controlled
by
the celebrated 
Denjoy-Carleman-Ahlfors theorem \cite{Hay7,Nev}. 

\begin{thm}[Denjoy-Carleman-Ahlfors theorem]
 \label{thmdca}
Let   $f$ be a transcendental meromorphic function in the plane of finite lower order $\mu$.
Then the number of direct transcendental singularities of $f^{-1}$  is at most $\max \{1, 2 \mu \}$.
\end{thm}

A key consequence of Theorem \ref{thmdca} is that a transcendental entire function of finite lower order $\mu$ has at most $2 \mu $ 
finite asymptotic values \cite{Hay7}. A  result of Bergweiler and Eremenko \cite{BE} shows that the critical values of a meromorphic 
function of finite (lower) order have a decisive influence on 
indirect transcendental singularities. 

\begin{thm}[\cite{BE}]
 \label{bethm}
Let $f$ be a transcendental meromorphic function in the plane of finite lower order.\\
(a) If $f^{-1}$ has an indirect transcendental singularity over $\alpha \in \widehat \C = \C \cup \{ \infty \}$
then  each   neighbourhood  of the singularity  contains infinitely many zeros of $f'$ which are not $\alpha$-points of $f$;
in particular, $\alpha$ is a limit point of critical values of $f$. \\
(b) If $f$ has finitely many critical values then 
$f^{-1}$ has finitely many transcendental singularities,  and all transcendental singularities are logarithmic.  
\end{thm}

Theorem \ref{bethm} was proved in \cite{BE} for $f$ of finite  order, and was extended to finite lower order, using essentially the same method,
by Hinchliffe \cite{Hin}.
Part (b) follows from  part (a) combined with Theorem \ref{thmdca} and a well known  classification theorem  from \cite[p.287]{Nev},
which shows in particular that any transcendental singularity of the inverse function over an isolated singular value is logarithmic.
Theorem \ref{bethm} was employed in \cite{BE} to prove a long-standing conjecture of Hayman \cite{Hay1}
concerning zeros of $ff'-1$, and has found many subsequent
applications, including to zeros of derivatives \cite{La11}. 
The reader is referred to \cite{BEdir,sixsmith} for further striking results on singularities of the inverse, both restricted to entire functions 
but independent of the order of growth.

The starting question of the present paper concerns the extent to which 
Theorems \ref{thmdca} and \ref{bethm} 
hold under the weaker hypothesis that
$f^{(k)}/f$ has finite lower order for some $k \in \N = \{ 1, 2, \ldots \}$. The obvious example $f(z) = \exp( \exp(z) )$ shows that $f'/f$ can have finite
order despite $f$ having infinite lower order; here $f^{-1}$ has infinitely many direct (indeed logarithmic) singularities over $0$ and $\infty$, and one over~$1$.
Furthermore, if $k \in \N$  and $A_k$ is a transcendental entire function then the lemma of the logarithmic derivative \cite{Hay2}
shows
that every
non-trivial solution of 
\begin{equation}
 \label{de1}
w^{(k)} - A_k(z) w = 0
\end{equation}
has infinite lower order, even if  $A_k$ has finite  order. Clearly each of $\exp( \exp(z) )$ and $\exp( z^{-1} \sin z )$ satisfies an equation of form 
(\ref{de1}) with coefficient of finite order. Note further that if $f$ is a transcendental meromorphic function in the plane 
and $f'/f$ has finite lower order then it is easy to prove by induction that so has $A_k = f^{(k)}/f$ for every $k \geq 1$, using the formula
$A_{k+1} = A_k A_1 + A_k'$, whereas the example 
$$
f(z) = e^{-z/2} \sin (e^z), \quad \frac{f'(z)}{f(z)} = - \, \frac12 + e^z \cot (e^z), \quad \frac{f''(z)}{f(z)} = \frac14 - e^{2z} ,
$$
shows that $f''/f$ can have finite order despite $f'/f$ having infinite lower order. 

\begin{thm}
 \label{thmds}
Let  $f$ be  a transcendental meromorphic function in the plane such that  $f^{-1}$  has $n \geq 1$ distinct direct transcendental singularities over  finite non-zero values. 
Let $k \in \N$ and let $\mu$ be the lower order of $A_k = f^{(k)}/f$.
Then the following statements hold.
\\
(i) There exists  a set $F_0 \subseteq [1, \infty)$ of 
finite logarithmic measure such that 
\begin{equation}
 \label{b5a}
\lim_{r \to + \infty, r \not \in F_0} \frac{ \log \left( \min \{ |A_k(z)| : |z| = r \} \right)}{\log r} = - \infty .
\end{equation}
(ii) If $n \geq 2$ then $n \leq 2 \mu $. \\
(iii) If $n=1$  and there exist $\kappa > 0$ 
and  a path $\gamma$ tending to infinity in 
the complement of the  neighbourhood $C(\kappa)$ of the singularity, then $\mu \geq 1/2$. 
\end{thm}
Theorem \ref{thmds} will be deduced from a version of the Wiman-Valiron theory for meromorphic functions with direct tracts developed in \cite{BRS},
and part (ii) is sharp, by Example I in Section \ref{exa}.
Furthermore, if $g$ is a transcendental entire function of lower order less than $ 1/2$
then the inverse function of  $f=1 - 1/g$ has a direct singularity over $1$; in this case $A_k$ obviously has lower order less than $ 1/2$
but the $\cos \pi \lambda$ theorem 
\cite[Ch. 6]{Hay7} implies that 
every neighbourhood of the singularity contains  circles $|z| = r$ with $r$ arbitrarily large, so that a path $\gamma$ as in (iii) cannot exist.

\begin{thm}
 \label{bekeylem2}
Let $f$ be a transcendental meromorphic function in the plane such that $f^{(k)}/f$ has finite lower order for some $k \in \N$. 
Assume that  $f^{-1}$ has an indirect transcendental singularity over $\alpha \in \widehat \C $. Then each
neighbourhood of the singularity  contains infinitely many zeros of $f'f^{(k)}$ which are not $\alpha$-points 
of $f$.  
\end{thm}
Theorem    \ref{bekeylem2} will be proved using a modification of methods from \cite{BE,Hin}.

\begin{cor}
 \label{cor1}
Let $f$ be a transcendental meromorphic function in the plane, with finitely many critical values, such that $f'/f$ has finite  lower order.
Then $f^{-1}$ has  finitely many transcendental singularities over finite non-zero values,  and $f$ has finitely many asymptotic
values. Moreover,  all transcendental singularities of $f^{-1}$ are logarithmic.
\end{cor}
Corollary \ref{cor1} follows from  Theorems \ref{thmds} and \ref{bekeylem2}, coupled with \cite[p.287]{Nev}. 

\begin{cor}
 \label{cor3}
Let $f$ be a transcendental meromorphic function in the plane such that $f''/f$ has lower order 
$\mu < \infty$ 
and $f'/f$ and $f''/f'$ have finitely many zeros.
Then $f''/f'$ is a rational function and $f$ has finite order and finitely many poles.
\end{cor}
To prove Corollary \ref{cor3} observe  that all but finitely many zeros of $f'f''$ are zeros of $f$.
Thus  $f^{-1}$ has no indirect  singularities, by Theorem \ref{bekeylem2}, 
and hence $f$ has finitely many asymptotic values, in view of Theorem \ref{thmds}. Since 
$f$ evidently has finitely many critical values, 
the result follows via
\cite[Theorem 2]{La11}. The condition $\mu < \infty$ holds if $f'/f$ has finite lower order, 
and is not redundant, because of an example in \cite{La11}.
\hfill$\Box$
\vspace{.1in}

The last result of this paper is related to the following theorem from \cite{Laschwarzian}. 

\begin{thm}[\cite{Laschwarzian}]\label{thm1}
Let $M$ be a positive integer and let $f$ be a transcendental meromorphic function in the plane
with transcendental Schwarzian derivative 
\begin{equation}
 \label{Schdef}
S_f(z) = \frac{f'''(z)}{f'(z)} - \frac32 \left( \frac{f''(z)}{f'(z)} \right)^2  ,
\end{equation}
such that:
(i) $f$ has finitely many critical values and
 all multiple points of $f$ have multiplicity at most $M $;
(ii) the inverse function of $f$ has finitely many transcendental singularities.

Then the following three conclusions hold:
(a)  $f$ has infinitely many multiple points;
(b) the inverse function of $S_f$ does not have a direct transcendental singularity over $\infty$;
(c) the value $\infty$ is not Borel exceptional for $S_f$.
 \end{thm}
Conclusion (a) is a result of Elfving \cite{elf} and Rolf Nevanlinna \cite{Nev2,Nev}, but was proved in \cite{Laschwarzian} by
a completely different method.
The following example shows that under the hypotheses of Theorem~\ref{thm1}
the inverse of the Schwarzian can have a direct transcendental
singularity over a finite value: write 
$$
g(z) = \sinh z , \quad S_g(z) = 1 - \frac{3 \tanh^2 z}2 ,
$$
so that  $S_g^{-1}$ has two  logarithmic singularities over $-1/2$. 
However, assumptions (i) and (ii) of  Theorem \ref{thm1} imply that  $f$ belongs to the Speiser class $\mathcal{S}$ \cite{Ber4,BE} 
consisting of all meromorphic functions in the plane for which the inverse function has finitely many singular values. 
For $f \in \mathcal{S}$, the following  result excludes direct singularities of the inverse of $S_f$ over $0$.

\begin{thm}\label{thm1aa}
Let  $f$ be a transcendental meromorphic function in the plane belonging to the Speiser class $\mathcal{S}$, with transcendental Schwarzian derivative $S_f$.
Then the inverse function of $S_f$  does not have a direct transcendental singularity over $0$.
\end{thm}
The example $f(z) = \tan^2 \sqrt{  z}$ from  \cite{CER} shows that for $f \in \mathcal{S}$ it is possible for $0$ to be an asymptotic value of $S_f$. 
Here direct computation shows that 
$f''(z)/f'(z)$ tends to $0$ as $z \to \infty$ in the left half plane, and  so does $S_f(z)$.

The author thanks the referees for helpful comments. 

\section{Examples illustrating Theorems \ref{thmds} and \ref{bekeylem2}}\label{exa}

\noindent
\textbf{Example I.} A function extremal for 
Theorem \ref{thmds}(ii), but not for Theorem \ref{thmdca}, is given by 
$$
f(0) = 1, \quad 
\frac{f'(z)}{f(z)} = \frac{\pi z}{\sin \pi z} .
$$
Here $f$ is meromorphic in the plane,
having at each non-zero integer $n$ a zero or pole of multiplicity $|n|$, depending on the sign and parity of $n$.
Hence  $N(r, f) $ and $N(r, 1/f)$ have order $2$. 
Because
$$
0 < 
\alpha = 
\int_0^{+\infty} \frac{\pi y}{\sinh \pi y} \, dy =
\int_0^{+\infty} \frac{\pi y}{\pi y + (\pi y)^3/6 + \ldots } \, dy < \int_0^1 1 \, dy + \int_1^\infty \frac{6}{\pi^2  y^2  } \, dy < \pi 
$$
and  $f'/f$ is even, $f$ has distinct  asymptotic values $e^{\pm i \alpha }$,
approached as $z$ tends to 
infinity along the imaginary axis.   As $f'/f$ has finite order and $f$ has no finite non-zero critical values, 
both  of these singularities of $f^{-1}$ are direct by Theorem \ref{bekeylem2}.
\hfill$\Box$
\vspace{.1in}
\\\\
\textbf{Example II.} Define $g$ by
$$
g(0) = 1, \quad 
\frac{g'(z)}{g(z)}  = A_1(z) =  \frac1{\pi \cos \sqrt{z}} .
$$
The zeros of $\cos \sqrt{z}$ occur where $\sqrt{z} = b_n = (2n+1) \pi/2 $, with $n \in \Z$, and the residue of $A_1$ at $ b_n^2$ is $\pm (2n+1)$. 
Thus $g$ is meromorphic in $\C$, with  zeros and poles in $\R^+$ and no finite non-zero critical values. 
Integration along the negative real axis shows that $g$ has a  non-zero real asymptotic value $\alpha$,
and $g^{-1}$ has  a logarithmic singularity over 
$\alpha$ by Corollary \ref{cor1}.  
This gives
$\delta > 0$  and a simply connected component $C$ of 
$\{ z \in \C: |g(z) - \alpha | < \delta \}$ with  $(-\infty, R) \subseteq C$ for some $R < 0$.
Moreover,  $C$ is symmetric with respect to
$\R$, since $g$ is real meromorphic, so that $C \cap \R^+$ is bounded, and $g$ is extremal for  Theorem \ref{thmds}(iii).
\hfill$\Box$
\vspace{.1in}

\noindent
\textbf{Example III.} Let $F(z) = \exp( -z /2 - (1/4) \sin 2z ) \cos z $, so that $F''/F$ is entire of finite order.
Then $F(z)$  tends to $0$ along $\R^+$ and this singularity of $F^{-1}$  is evidently indirect. 
\hfill$\Box$
\vspace{.1in}

\noindent
\textbf{Example IV.} Define entire functions $A_1$ and $v$ by
$$
v(0) = 1, \quad 
\frac{v'(z)}{v(z)} = A_1(z) = \frac{1 - \cos z}{z^2} = \frac12 + \ldots . 
$$
Then there exists $\alpha \in \R^+$ such that $v(x) \to \exp( \pm \alpha )$ as $x \to \pm \infty$ on $\R$
and, since $A_1$ does not satisfy (\ref{b5a}), 
Theorem \ref{thmds} 
implies that $v^{-1}$ has no direct singularities over finite non-zero values.
Because all critical points of $v$ are real, all but finitely many of them belong to neighbourhoods of 
the indirect singularities over $\exp( \pm \alpha )$, and so $v^{-1}$ has no other indirect singularities, by Theorem~\ref{bekeylem2}.
Thus applying \cite[p.287]{Nev} again, in conjunction with Iversen's theorem,
shows that $v^{-1}$ has logarithmic singularities over the omitted values $0$ and $\infty$. 
\hfill$\Box$
\vspace{.1in}

\noindent
\textbf{Example V.} Let $h(z) = \exp(  \sin z - z )$, so that $A_1 = h'/h$ is
entire of finite order but does not satisfy (\ref{b5a}). Since $h(z)$ tends to $0$ along  $\R^+$, and to $\infty$ on the negative real axis,
with $h'(2 \pi n )= 0$ for all $n \in \Z$, these 
singularities of $h^{-1}$  are direct but not logarithmic.
\hfill$\Box$
\vspace{.1in}


\section{Preliminaries}

The following well known estimate may be found in  Theorem 8.9 of \cite{Hay7}.

\begin{lem}[\cite{Hay7}]
 \label{dcalem}
Let 
$D_1, \ldots , D_n$ be $n \geq 2$ pairwise disjoint plane
domains.
If $u_1, \ldots, u_n $ are non-constant subharmonic functions on $\C$ such that
$u_j$ vanishes outside $D_j$  then 
\begin{equation}
 \label{dca0}
\liminf_{r \to \infty} \frac{h(r)}{r^{n/2} } > 0, \quad h(r) = \max_{1 \leq j \leq n}  B(r, u_j), \quad B(r, u_j) = \sup \{ u_j(z): |z| = r \}  .
\end{equation}
\end{lem}
\hfill$\Box$
\vspace{.1in}

For $a \in \C$ and $R > 0$ 
denote by $D(a, R)$ the open disc of centre $a$ and radius $R$, and by $S(a, R)$ its boundary circle. 

\begin{lem}
 \label{lemderivs}
To each 
$k \in \N$ corresponds $d_k \in (0, \infty) $ with the following property. 
Suppose that $0 < R < \infty $ and $w=h(z)$ maps the domain $U \subseteq \C$ conformally onto  $D(a, R)$, with inverse function $F: D(a, R) \to U$. Then 
there exists an analytic function $V_k: D(a, R) \to \C$ with 
\begin{equation}
 \label{derivest1}
h^{(k)}(z) F'(w)^k = V_k(w), \quad 
|V_k(w)| \leq \frac{d_k}{( R-|w-a|)^{k-1} }  \quad \hbox{as} \quad |w-a| \to R- .
\end{equation}
\end{lem}
\textit{Proof.} Assume that $a=0$ and initially that $R=1$. It 
is  clear that (\ref{derivest1}) holds for $k=1$, with $V_1 (w) = 1$. If  (\ref{derivest1}) holds
for $k$ then it follows that 
\begin{equation*}
h^{(k+1)}(z) F'(w)^{k+1} = V_k'(w) - k h^{(k)}(z) F'(w)^{k-1} F''(w) = V_k'(w) - k V_k(w) \, \frac{ F''(w) }{F'(w)} .
\end{equation*}
Since $F''(w)/F'(w) = O( 1-|w| )^{-1}$ as $|w| \to 1-$ by \cite[p.5, (1.6)]{Hay9}, applying Cauchy's estimate for derivatives to $V_k$ 
proves the lemma by induction when $R=1$. 
In the general case write $w = h(z) = RH(z) = Rv$ and $z = F(w)= G(v)$ so that, as $|w| \to R-$,  
$$|h^{(k)}(z) F'(w)^k | = R^{1-k} |H^{(k)}(z) G'(v)^k |  \leq   \frac{d_k R^{1-k} }{( 1-|v|)^{k-1} }  = \frac{d_k}{( R-|w|)^{k-1} } .$$
\hfill$\Box$
\vspace{.1in}

\begin{lem}
 \label{belem6a}
Let $M \in \N$ and $s > 2^{24}$ and let $E_1, \ldots , E_{N}$
be $N \geq 24M$ pairwise disjoint domains in $\C$,
and for $t > 0$ let $\phi_j(t)$ be the angular measure of 
$S(0, t) \cap E_j$. Then 
at least $N-12M$ of the 
$E_j$ satisfy
\begin{equation}
\int_{ [ 4s^{1/2} , s /4 ] } \frac{\pi \, dt}{t \phi_j(t)} > M \log s 
\quad \hbox{and} \quad 
\int_{ [ 4s , s^2 /4 ] } \frac{\pi \, dt}{t \phi_j(t)} > M \log s .
\label{Djnarrow}
\end{equation}
\end{lem}
\textit{Proof.} 
This is a standard application as in \cite[Ch. 8]{Hay7} or \cite{BE} of 
the Cauchy-Schwarz inequality, which gives
\begin{equation}
 \label{csineq}
\frac{L^2}{t} \leq \frac1t
\left( \sum_{j=1}^{L} \phi_j(t) \right)
\left( \sum_{j=1}^{L} \frac1{  \phi_j(t) }\right)
\leq 2  \sum_{j=1}^{L} \frac{\pi}{ t \phi_j(t) }
\end{equation}
for $M \leq L \leq N$ and $t > 0$.
If $s > 2^{24}$
and either inequality of (\ref{Djnarrow}) fails for $L \geq 6M $ of the $E_j$, without loss of generality for
$j=1, \ldots, L$,  then integrating (\ref{csineq}) 
yields a contradiction via 
$$
2 L M \log s < 6 L M \log \frac{\sqrt{s}}{16} \leq L^2 \, \log \frac{\sqrt{s}}{16} \leq 2L M \log s  .
$$
\hfill$\Box$
\vspace{.1in}

\begin{lem}[\cite{Ber4}]
\label{lember4}
Let $h$ be a transcendental meromorphic function in the plane belonging to the Speiser class $\mathcal{S}$. 
Then there exist positive constants $C$, $R$ and $M$ such that 
\begin{equation}
 \label{Sineq}
\left| \frac{z h'(z)}{h(z)} \right| \geq C \log^+ \left| \frac{h(z)}M \right| \quad \text{for $|z| \geq R$.}
\end{equation}
 \end{lem}

\section{Proof of Theorem \ref{thmds}} 

Let   $f$ be a transcendental meromorphic function in the plane such that $f^{-1}$ has $n \geq 1$ 
direct singularities  over  (not necessarily distinct) finite non-zero values $a_1, \ldots, a_n$. Let $k \in \N$; then $A_k = f^{(k)}/f$ does not vanish identically.
There exist a small positive $\delta $ and 
non-empty components
$D_j $ of  $\{ z \in \C : |f(z) - a_j| < \delta \}$, for  $j = 1, \ldots n$,
such that $f(z) \not = a_j$ on $D_j$,  so that $D_j$ immediately
qualifies as a direct tract for $g_j = \delta /(f-a_j)$ in the sense of \cite[Section 2]{BRS}. 
Here $\delta$ may be chosen so small that if $n \geq 2$ then these $D_j$ are pairwise disjoint.
For each $j$,  define a non-constant subharmonic function $u_j$ on $\C$ by
$$
u_j(z) = \log |g_j(z)| = \log \left| \frac{\delta }{ f(z) - a_j} \right|  \quad 
(z \in D_j ), \quad
u_j(z) = 0  \quad ( z \not \in D_j).
$$
Then  \cite[Theorem 2.1]{BRS} implies that, with $B(r, u_j)$ as in (\ref{dca0}), 
\begin{equation}
 \label{uest1aa}
\lim_{r \to + \infty} \frac{B(r, u_j)}{\log r} = + \infty, \quad 
\lim_{r \to + \infty} a(r, u_j) = + \infty, \quad a(r, u_j) = rB'(r,u_j) .
\end{equation}


\begin{lem}
 \label{lemwv}
There exists a set  $F_0 \subseteq [1, \infty)$, of finite logarithmic measure, such that for each $s \in [1, \infty) \setminus F_0$ 
and each $j$ there exists $z_j$ with 
\begin{equation}
\label{wvest}
|z_j| = s, \quad A_k(z_j) = 
 \frac{f^{(k)}(z_j)}{f(z_j)} =  O \left(  \exp( - B(s, u_j) /2 )   \right).
\end{equation}
\end{lem}
\textit{Proof.}
Fix $\tau$ with $1/2 < \tau < 1$ and apply the version of Wiman-Valiron theory
developed in \cite{BRS} for meromorphic functions with
direct tracts. 
By \cite[Theorem 2.2 and Lemma 6.10]{BRS}  there exists a set $F_0 \subseteq [1, \infty)$ of finite logarithmic measure such that every
$s  \in [1, \infty) \setminus F_0$ has the following two properties: first, $a(s, u_j) $ is large, by (\ref{uest1aa}), but satisfies
\begin{equation}
 \label{dca5}
a(s, u_j) \leq B(s, u_j)^2 ; 
\end{equation}
second, for each $j$ there exists $z_j $ with $|z_j| = s$ and $u(z_j) = B(s, u_j)$ such that 
\begin{equation}
 \label{dca6}
\frac{f(z)-a_j}{f(z_j)-a_j} \sim \left( \frac{z}{z_j} \right)^{-a(s, u_j)} 
\quad \hbox{for} \quad |z-z_j| < \frac{s}{ a(s,u_j)^{\tau} } .
\end{equation}
A standard application of Cauchy's estimate for derivatives in (\ref{dca6}) now gives 
\begin{equation*}
\left( \frac{f'}{f-a_j} \right)^{(p)} (z) = 
O \left(  \frac{ a(s, u_j) }{s} \right)^{p+1}
\quad \hbox{for} \quad p=0, \ldots, k-1 \quad \hbox{and} \quad |z-z_j| < \frac{s }{2 a(s,u_j)^{\tau} } .
\end{equation*}
It follows via \cite[Lemma 3.5]{Hay2} that
$$
\frac{f^{(k)}(z_j)}{f(z_j)} =  \frac{f^{(k)}(z_j)}{f(z_j) - a_j}  \cdot \frac{f(z_j)-a_j}{f(z_j)}  
=  O \left(  \frac{ a(s, u_j)^{k} \exp( - B(s, u_j) )  }{s^k}  \right),
$$
which, by (\ref{dca5}), yields  (\ref{wvest}) for large enough $s \not \in F_0$.
\hfill$\Box$
\vspace{.1in}

Combining (\ref{uest1aa}) with  (\ref{wvest}) for $j=1$ leads to (\ref{b5a}). 
To prove the remaining assertions it may be assumed that $A_k$ has finite lower order $\mu$. 
Choose a positive sequence $(r_m)$ tending to infinity such that 
\begin{equation}
 \label{dca2}
T(8r_m, A_k) < r_m^{\mu+o(1)} .
\end{equation}
Let $m$ be large and let $w_1, \ldots, w_{q_m}$ be the zeros and poles of $A_k$ in $r_m/4 \leq |z| \leq 4r_m$, repeated according to multiplicity:
then (\ref{dca2}) and standard estimates yield
\begin{equation}
 \label{dca3}
q_m \leq n(4r_m, A_k) + n(4r_m,1/A_k) \leq \frac2 {\log 2} \,  T(8r_m, A_k) + O(1) \leq r_m^{\mu+o(1)} . 
\end{equation}
Let $U_m$ be the union of the discs $D(w_j, r_m^{ - \mu })$. Since the sum of the radii of the discs of $U_m$ is  $ o( r_m)$ by (\ref{dca3}),  there 
exists a set $E_m \subseteq [r_m/2, 2r_m ]$, of linear measure at least $r_m$, and so logarithmic measure $l_m \geq 1/2$,
such that for $r \in E_m$ the circle $|z| = r$ does not meet $U_m$. 
A standard application of the Poisson-Jensen formula \cite{Hay2} on the disc 
$|\zeta| \leq 4r_m$ 
then yields 
\begin{equation}
 \label{dca4}
\left| \log \left| A_k(z) \right| \right| \leq r_m^{\mu+o(1)} \quad \hbox{for} \quad |z| \in E_m .
\end{equation}
Since $m$ is large and $l_m \geq 1/2$,
there exists $s_m \in E_m \setminus F_0$.

Suppose now that   $n=1$ and  there exist $\kappa > 0$ 
and  a path $\gamma$ tending to infinity in 
the complement of the  neighbourhood $C(\kappa)$ of the singularity, or that $n \geq 2$.
Then (\ref{dca0}) holds, by  \cite[Theorem 6.4]{Hay7} when $n=1$, and by Lemma \ref{dcalem} when $n \geq 2$.
Combining (\ref{dca0}) and (\ref{wvest}), with $s= s_m \geq r_m /2$, 
yields points $z_j$ with $|z_j| = s_m$ and, for at least one $j$, 
$$A_k(z_j) = O \left(  \exp( - B(s_m, u_j) /2 )   \right) = O \left( \exp \left( -s_m^{n/2-o(1)} \right) \right) .$$
On combination with (\ref{dca4}) this forces $2 \mu \geq n$. 
\hfill$\Box$
\vspace{.1in}

\section{Indirect singularities} \label{transing}

\begin{prop}
 \label{bekeylem}
Let $f$ be a transcendental meromorphic function in the plane such that $ f^{(k)}/f$ has finite lower order $\mu$ for some $k \in \N$. 
Assume that $f^{-1}$ has an indirect transcendental singularity over $\alpha \in \C \setminus \{ 0 \}$. 
Then for each $\delta > 0$ the  neighbourhood $C(\delta)$ of the singularity  contains infinitely many zeros of $f' f^{(k)}$.  
\end{prop}

The proof of Proposition \ref{bekeylem} will take up the whole
of this section. 
The method is adapted from those in \cite{BE,Hin}, but some complications arise, in particular when $k \geq 2$. 
Assume throughout that $f$ and $\alpha$ are as in the hypotheses but 
$C( \varepsilon )$, for some small $\varepsilon > 0$, contains
finitely many zeros of $f'f^{(k)}$. It may be assumed that $\alpha =1$, and that 
$C( \varepsilon )$ contains no zeros of $f'f^{(k)}$. 
Choose positive integers $N_1, N_2, \ldots, N_9 $ with $5 \mu  + 12 < N_1$ and $N_{j+1}/N_j$ large for each $j$.

\begin{lem}
 \label{belem2}
For each $j \in \{  1, \ldots, N_9 \}$ there exist  $z_j \in C( \varepsilon )$
and  $a_j \in \C$ with
$0 < r_j = | 1 - a_j|  < \varepsilon /2 $, as well as a simply connected domain $D_j \subseteq C( \varepsilon )$,
with the following properties. 
The $a_j$ are pairwise distinct and the $D_j$ pairwise 
disjoint.
Furthermore, the function $f$ maps $D_j$ univalently onto $D(1, r_j)$, with $z_j \in D_j$ and 
$f(z_j) = 1$. Moreover,  $0 \not \in D_j$ but $D_j$ contains a path $\sigma_j $  tending to infinity, which
is mapped by
$f$ onto the half-open line segment $[1, a_j)$, with $f(z) \to a_j$ as $z \to \infty $ on $\sigma_j$.
\end{lem}
This is proved exactly as in \cite{BE}. If $0 < T_j < \varepsilon /2$ and $z_j \in C(T_j)$ is such that 
$f(z_j) = 1$, let $r_j$ be the supremum of $t > 0$ such that the branch of $f^{-1}$ mapping $1$ to $z_j$ admits unrestricted analytic continuation in 
$D(1, t)$. Then $r_j < T_j$ because $f$ is not univalent on $C(T_j)$, and there is a singularity $a_j$ of $f^{-1}$ with $|1-a_j| = r_j$; moreover,
$a_j$ must be an asymptotic value of $f$. The $z_j$ and $T_j$ are then chosen inductively: for the details see \cite{BE}
(or \cite[Lemma 10.3]{lajda}).
\hfill$\Box$
\vspace{.1in}

\begin{lem}\label{belem3}
 Let the $z_j, a_j, \sigma_j$ and
$D_j$ 
be as in Lemma $\ref{belem2}$.
For $t > 0$, let $t \theta_j(t) $ be the length of the longest open arc
of $S(0, t)$ which lies in $D_j$. Then $f$ satisfies, as $z$ tends to infinity on $\sigma_j$,
\begin{equation}
\log \frac{r_j}{|f(z) - a_j|} \geq 
\int_{|z_j|}^{|z|} \frac{dt}{4 t \theta_j(t) }  .
\label{be1}
\end{equation}
\end{lem}
\textit{Proof.} 
Let $z = H(w)$ be the branch of $f^{-1}$ mapping
$D(1, r_j )$ onto $D_j$. For $z \in \sigma_j$ the distance from $z$ to
$\partial D_j$ is at most $|z| \theta_j(|z|)$. Thus  Koebe's quarter theorem \cite[Ch. 1]{Hay9} 
implies that
$$
|(w - a_j ) H'(w)| \leq 4 |z| \theta_j(|z|) 
\quad \hbox{for} \quad z = H(w), \, w \in [1 , a_j) .
$$
Hence, for large $z \in \sigma_j$ and $w = f(z)$,
writing $u = H(v)$ for $v \in [1, w]$ gives (\ref{be1}) via
\begin{eqnarray*}
\log  \frac{r_j}{| f(z)-a_j |}  
&=& \int_{[1, w]} \frac{|dv|}{|a_j-v|}
= \int_{H([1,w])} \frac{|du|}{|(a_j-v)H'(v)|} \geq \int_{H([1,w])} \frac{|du|}{4 |u| \theta_j(|u|)} .
\end{eqnarray*} 
\hfill$\Box$
\vspace{.1in}
 
Since $N_1 > 5 \mu $ there exists  a positive sequence $(s_n)$ tending to infinity
such that
\begin{equation}
T(s_n^5, f^{(k)}/f) + T(s_n^5 , f/f^{(k)}) \leq s_n^{N_1}  .
\label{Mexist}
\end{equation}
Set
\begin{equation}
G(z) = z^N \frac{f^{(k)}(z)}{f(z)}, \quad N = N_5. 
\label{Gdefn}
\end{equation}
Applying \cite[Lemma 4.1]{LaRo2} to $1/G$ (with 
$\psi (t) = t$ in the notation of \cite{LaRo2}) gives a small positive $\eta$ such that 
$G$ has no critical values $w$ with $|w| = \eta $ and such that
the length $L(r, \eta , G)$ of the level curves $|G(z)| = \eta $ lying in $D(0, r)$ satisfies
\begin{equation}
L(s_n^4, \eta, G) = O( s_n^6 T( e^8  s_n^4, G)^{1/2} ) =
O( s_n^{6 + N_1/2} ) \leq s_n^{N_1}  \quad \hbox{ as $n \to \infty$,}
\label{lrgest}
\end{equation}
using (\ref{Mexist}) and the fact that $N_1 > 12$. 
Assume henceforth that $n$ is large.

\begin{lem}
 \label{lemingrowth}
At least $N_8$ of the domain $D_j$ and paths $\sigma_j$, without loss of generality 
$D_1, \ldots, D_{N_8}$ and 
$\sigma_1, \ldots, \sigma_{N_8}$, are such that 
\begin{equation}
 \label{minest1}
|f(z)-a_j| \leq  s_n^{- N_7} \quad \hbox{for $z \in \sigma_j$ with $|z| \geq s_n/4$.}
\end{equation}
\end{lem}
\textit{Proof.} 
By Lemma \ref{belem6a} it may be assumed that, for $j = 1, \ldots, N_8$,  
\begin{equation*}
 \int_{ [ 4s_n^{1/2} , s_n /4 ] } \frac{\pi \, dt}{t \theta_j(t)} > N_8 \log s_n ,
\end{equation*}
which, on combination with Lemma \ref{belem3}, leads to
 (\ref{minest1}). 
\hfill$\Box$
\vspace{.1in}

\begin{lem}
 \label{belem5a}
Let $w_1, \ldots , w_{q_n}$ be the zeros and poles of $G$ in $s_n^{1/4} \leq |z| \leq s_n^4$,
repeated according to multiplicity. Then 
\begin{equation}
 \label{n(r)est1}
q_n \leq  n(s_n^4, 1/G) + n(s_n^4, G) = o \left( s_n^{N_1} \right) 
\end{equation}
and there exist $t_n , T_n$ satisfying
\begin{equation}
s_n^{1/2} - 1 \leq t_n \leq s_n^{1/2} , \quad
s_n^2 \leq T_n \leq s_n^2 +1, 
\label{tnTndef}
\end{equation}
such that
\begin{equation} 
\max \{ | \log |G(z)| | : z \in S(0, t_n) \cup S(0, T_n) \} 
\leq s_n^{N_1+1} .
\label{tnTn2}
\end{equation}
\end{lem}
\textit{Proof.} (\ref{n(r)est1}) follows from  (\ref{Mexist}).
Let $U_n$ be the union of the discs $D(w_q, s_n^{-N_1-1} )$: these discs  have sum of radii at
most $s_n^{-1}$ and so since $n$ is large there exist $t_n, T_n$
satisfying (\ref{tnTndef})  such that
the circles $S(0, t_n), S(0, T_n)$ do not meet $U_n$. Hence the Poisson-Jensen formula gives
(\ref{tnTn2}).
\hfill$\Box$
\vspace{.1in}

\begin{lem}
\label{lemnocomp}
Define sets $E$, $K_n$ and $L_n$ by $E = \{ z \in \C : |G(z)| < \eta \} $ and 
$$
K_n = \{ z \in \C : t_n < |z| < T_n \} , \quad L_n = \{ z \in \C : s_n/4 < |z| < 4 s_n \}.
$$
Then the number of components $E_q$ of $E \cap K_n $ which meet 
$L_n $ is at most $s_n^{N_1}$.
\end{lem}
\textit{Proof.} If the closure $F_q$ of $E_q$ lies in $K_n $ then 
$E_q$ must contain a zero of $G$, whereas if $F_q \not \subseteq K_n $ then $\partial E_q \cap K_n$ 
has arc length at least $s_n/8$. Thus the lemma follows from (\ref{lrgest}) and (\ref{n(r)est1}). 
\hfill$\Box$
\vspace{.1in}

\begin{lem}
 \label{belem5}
Let $u$ lie on $\sigma_j$ with $s_n/4 \leq |u| \leq 4s_n$.
Then, with $d_k$ as in Lemma \ref{lemderivs}, there exists 
$v$ on $\sigma_j$
such that:
\begin{equation}
 \label{vest1}
|u| \leq |v| \leq  |u| + s_n^{-N_3}; \quad 
|f(v ) - a_j| \leq | f(u) - a_j | ;  \quad 
|f^{(k)}(v) | \leq k^k  d_k s_n^{kN_3} |f(u)-a_j| .
\end{equation}
\end{lem}
\textit{Proof.} Starting at $u$, follow $\sigma_j$ in the direction
in which $|f(z) - a_j|$ decreases. Then $\sigma_j$ describes an arc 
$\gamma$ joining the circles $S(0, |u|)$ and $S(0, |u| + s_n^{-N_3})$, such that
the first two inequalities of (\ref{vest1}) hold for all 
$v \in \gamma$. 
Since $f$ maps $D_j$ univalently onto $D(1, r_j)$, the inverse function $H$ of $f$ maps
a proper sub-segment $I$ of the half-open line segment $J = [ f(u)  , a_j  ) $ onto $\gamma$. Assume that 
the last inequality of (\ref{vest1}) fails for all 
$v \in \gamma$. 
Then Lemma \ref{lemderivs} yields, on $I$, 
$$
|H'(w)| \leq k^{-1} s_n^{-N_3}|f(u)-a_j|^{-1/k}   (r_j-|w-1|)^{1/k-1}  .
$$
Since $1$, $f(u)$ and $a_j$ are collinear, a contradiction arises via
\begin{eqnarray*}
 s_n^{-N_3}   \leq \left| \int_I H'( w ) d w \right|  &\leq &  \int_I k^{-1} s_n^{-N_3} |f(u)-a_j|^{-1/k}   (r_j-|w-1|)^{1/k-1} \, |dw| \\
&< & \int_J k^{-1} s_n^{-N_3} |f(u)-a_j|^{-1/k}   (r_j-|w-1|)^{1/k-1} \, |dw| \\
&= & \int_{|f(u)-1|}^{r_j}  k^{-1} s_n^{-N_3} |f(u)-a_j|^{-1/k}   (r_j-t)^{1/k-1} \, dt \\
&=&  s_n^{-N_3} |f(u)-a_j|^{-1/k}   (r_j-|f(u)-1|)^{1/k} = s_n^{-N_3} .
\end{eqnarray*}
\hfill$\Box$
\vspace{.1in}

\begin{lem}
 \label{lemgronwall}
Let $E_p$ be a component of $E \cap K_n $ which meets $L_n$, and suppose that there exists $j = j(p)$ such that $E_p$  contains $k$ points $\zeta_1, \ldots , \zeta_k \in D_j$ 
each with 
$ |f(\zeta_q) - a_j| \leq  s_n^{-N_7} $. Assume further that $|\zeta_q - \zeta_{q'}| \geq s_n^{-N_3}$ for $q \neq q'$. 
Then $|f(z) - a_j | \leq s_n^{- N_2} $  for all $z \in E_p$, and $E_p \subseteq C(\varepsilon)$.
\end{lem}
\textit{Proof.} 
Let $M_0 = \sup \{ |f(z)| : z \in E_p \}$; then $M_0 < + \infty$ since poles of $f$ in $\C \setminus \{ 0 \}$ 
are poles of $G$, by (\ref{Gdefn}), and $|G(z)| \leq \eta $ on the closure of $E_p$. 
Choose $u_0 \in E_p$ with $|f(u_0) | \geq M_0/2$.
There exists a polynomial $P$, of degree at most $k-1$, such that 
$$
f(z) = P(z) + \int_{u_0}^z \frac{(z-t)^{k-1}}{(k-1)!}  f^{(k)}(t) \, dt \quad \hbox{on $E_p$. }
$$
The length of the boundary of $E_p$ is at most $2 s_n^{N_1}$ by (\ref{lrgest}). Hence 
each $z \in E_p$ can be joined to $u_0$ by a path in the closure of $E_p$, of length at most $4s_n^{N_1}$, and so 
\begin{equation}
 \label{flikeP}
|f(z) - P(z)| \leq M_0 \eta t_n^{-N_5} (2T_n)^{k-1} 4 s_n^{N_1} \leq M_0 s_n^{-N_4} ,
\end{equation}
by (\ref{Gdefn}) and (\ref{tnTndef}). 
In particular this gives $| P(\zeta_q) - a_j | \leq  (1+M_0) s_n^{-N_4} $
for each $q$. For  $z$ in $E_p$, Lagrange's interpolation formula leads to
\begin{eqnarray}
 \label{lagrange}
|P(z) - a_j| &=& \left| \sum_{q=1}^k (P(\zeta_q)-a_j) \prod_{\nu  \neq q} \frac{z-\zeta_\nu}{\zeta_q-\zeta_\nu} \right| \nonumber \\
&\leq & 
k  (1+M_0) s_n^{-N_4} (2T_n)^{k-1} s_n^{(k-1)N_3} \leq (1+M_0) s_n^{-N_3} .
\end{eqnarray}
Setting $z = u_0$ in (\ref{lagrange}) then delivers $M_0 \leq 2 |P(u_0)| \leq 2 |a_j| + o(1 + M_0) $ and so $M_0 \leq 5$. Now combining (\ref{flikeP}) with
(\ref{lagrange}) yields  $|f(z) - a_j | \leq s_n^{- N_2} $ and hence $|f(z) - 1| < \varepsilon$
on $E_p$. Since $E_p$ meets $D_j \subseteq C(\varepsilon)$, this gives  $E_p \subseteq C(\varepsilon)$. 
\hfill$\Box$
\vspace{.1in}

For each $j \in \{ 1, \ldots, N_8 \}$ choose $\lambda = s_n^{N_2}$ points $u_{j,1}, \ldots , u_{j,\lambda}$ on $\sigma_j$, each with 
$s_n/2 \leq |u_{j,\kappa} | \leq  s_n$ and such that $|u_{j,\kappa+1}| \geq |u_{j,\kappa} | + 2 s_n^{-N_3}$. 
Applying Lemma \ref{belem5} with $u = u_{j,\kappa} $ gives points $v_{j,\kappa} \in \sigma_j $ with 
$s_n/2 \leq |u_{j,\kappa}| \leq |v_{j,\kappa}| \leq |u_{j,\kappa} | +  s_n^{-N_3} \leq 2s_n$ and, using (\ref{Gdefn}), (\ref{minest1}) and (\ref{vest1}), 
\begin{equation}
 \label{vjest1}
|f(v_{j,\kappa}) - a_j|\leq s_n^{-N_7}, \quad 
|G(v_{j,\kappa})| \leq 2 |v_{j,\kappa} |^{N_5} | f^{(k)} (v_{j,\kappa})| \leq s_n^{-N_6} < \eta .
\end{equation}
These points $v_{j,\kappa}$ satisfy $|v_{j,\kappa+1}| \geq |v_{j,\kappa} | +  s_n^{-N_3}$, and each lies in a component of $E \cap K_n $ which meets $L_n$.
Since there are  $s_n^{N_2}$ of these  $v_{j,\kappa}$ for each $j$, but at most $s_n^{N_1}$ available components $E_p$ by Lemma \ref{lemnocomp},
it must be the case that for each $j$ there are at least $k$ points $v_{j,\kappa}$ lying in the same component $E_p$.
 Lemma \ref{lemgronwall} then implies that $E_p \subseteq C(\varepsilon)$ and $f(z) = a_j + o(1)$ on $E_p$.

Thus for $j=1, \ldots, N_8$ the following exist: a component $C_j = E_{p_j} \subseteq C(\varepsilon) $ of $E \cap K_n $ which meets $L_n$
and on which $f(z) = a_j + o(1)$; a point $v_j \in C_j$ such that, by  (\ref{vjest1}),
\begin{equation}
 \label{vjest2}
s_n/2 \leq |v_j|  \leq 2s_n , 
\quad 
|G(v_{j})| \leq   s_n^{-N_6} .
\end{equation}
Since $C_j \subseteq C(\varepsilon)$, the function $\log |1/G(z)|$ is subharmonic on $C_j$. Moreover, because $j' \neq j$ gives
$f(z) \to a_{j'} \neq a_j $ as $z \to \infty $ on $\sigma_{j'}$, 
the $C_j$ are pairwise disjoint and none of them contains a circle $S(0, t)$ with $t \in [t_n, T_n]$. 
For $t > 0$ let $\phi_j(t)$ be the angular measure of $C_j \cap S(0, t)$. 
Then  (\ref{tnTndef}) and  \cite[p.116]{Tsuji1}
give a harmonic measure  estimate   
$$
\omega (v_j, C_j, S(0, T_n) \cup S(0, t_n) ) \leq 
c_1 \exp \left( - \pi \int_{2|v_j|}^{T_n/2} \frac{dt}{t \phi_j(t)} \right) +
c_1 \exp \left( - \pi \int_{2t_n}^{|v_j|/2} \frac{dt}{t \phi_j(t)} \right) ,
$$
for $j=1, \ldots, N_8$, 
in which $c_1$ is a positive absolute constant.
By Lemma \ref{belem6a} and (\ref{tnTndef}), there exists at least one $j$ for which $ \omega (v_j, C_j, S(0, T_n) \cup S(0, t_n) ) \leq 2c_1 s_n^{-N_7}$.
For this choice of $j$ the two constants theorem  \cite{Nev} delivers,
using 
(\ref{tnTn2}), (\ref{vjest2})
and the fact that $|G(z)| = \eta $ on $\partial C_j \cap K_n$, 
$$
N_6 \log s_n \leq \log \frac1{|G(v_j)|} \leq \log \frac1\eta + 2c_1 
s_n^{- N_7 + N_1 + 1 }   ,
$$
a contradiction since  $n$ is large.
\hfill$\Box$
\vspace{.1in}

\section{Proof of Theorem \ref{bekeylem2}}

This is almost identical to the corresponding proof in \cite{BE}, but with Theorem \ref{thmds} standing in for the Denjoy-Carleman-Ahlfors 
theorem. 
Suppose that $f$, $k$ and $\alpha$ are as in the hypotheses but 
there exists 
$\varepsilon > 0$ such that in the neighbourhood $C( \varepsilon )$ of the singularity the function $f' f^{(k)}$ has finitely many  zeros 
which are not $\alpha$-points  of $f$: it may be assumed that there are no such zeros. 
On the other hand, because the singularity is indirect, $f$ must have infinitely many $\alpha$-points   in $C(\varepsilon)$. 
Since $f^{(k)}/f$ has finite lower order,
$f^{-1}$ cannot have infinitely many direct transcendental singularities over finite non-zero values,
by Theorem~\ref{thmds}. Set $A(\varepsilon) = \{ w \in \C : 0 <  | w- \alpha |  <  \varepsilon  \}$ if $\alpha \in \C$,
with $A(\varepsilon) = \{ w \in \C : |w| > 1/  \varepsilon  \}$ if $\alpha = \infty $.
In either case it may be assumed that $\varepsilon $ is so small that
$A(\varepsilon)  \subseteq \C \setminus \{ 0 \}$
and there is no $w$ in $A(\varepsilon)$
such that $f^{-1}$ has
a direct transcendental singularity over~$w$.

Take $z_0 \in C( \varepsilon )$,
with $f(z_0) = w_0 \neq \alpha $, and let  $g$ be that branch of
$f^{-1}$ mapping $w_0 $ to $z_0$. If $g$ admits unrestricted analytic continuation in $A(\varepsilon)$ 
then, exactly as in \cite{BE},  the classification theorem from   \cite[p.287]{Nev} shows that $z_0$ lies in a 
component $C_0$ of the set $\{ z \in \C :  f(z) \in A( \varepsilon ) \cup \{ \alpha \}  \}$ which 
contains at most
one point $z$ with $f(z) = \alpha$, so that 
$ C(\varepsilon) \not \subseteq C_0$. But any $z_1 \in C( \varepsilon )$ can be joined to $z_0$ by a path $\lambda$ on which $f(z) \in A(\varepsilon)\cup \{ \alpha \}$, 
which gives $\lambda \subseteq C_0$ and hence $ C(\varepsilon) \subseteq C_0$,  a contradiction.

Hence there exists a path $\gamma : [0, 1] \to A(\varepsilon)$,
starting at
$w_0$, such that 
analytic continuation of $g$ along $\gamma$ is not possible. 
This gives rise to
$S \in [0, 1]$ such that,
as $t \to S-$, the image 
$z = g( \gamma (t) )$ either
tends to infinity or to a zero $z_2 \in C( \varepsilon) $ of $f'$
with $f(z_2) = \gamma(S) \in A(\varepsilon )$, the latter  impossible by assumption. 
It follows that setting 
$z = \sigma(t) = g( \gamma (t) )$, for  $ 0 \leq t < S$, defines a path $\sigma$ 
tending to infinity
in $C( \varepsilon )$, on which $f(z) \to w_1 \in A(\varepsilon)$ as $z \to \infty$.
But then there exists $\delta > 0$ such that an unbounded subpath of $\sigma $ lies in a component
$C' \subseteq C( \varepsilon)$ of the set $\{ z \in \C : |  f(z) - w_1 | < \delta  \}$, with $\delta $ so small that  $f'f^{(k)}$ has no zeros
on $C'$. Further, the singularity over $w_1$ must be indirect,
since  direct singularities over values in $A(\varepsilon)$
have been excluded,
and this contradicts Proposition~\ref{bekeylem}.

\hfill$\Box$
\vspace{.1in}

\section{A result needed for Theorem \ref{thm1aa}}

\begin{thm}
[ \cite{LRW}, Theorem 1]
\label{thmlrw}
Let $u$ be  a subharmonic function in the plane 
such that $B(r)  = \sup \{ u(z): |z| = r \} $ satisfies 
$\lim_{r \to \infty} ( \log r)^{-1} B(r) = + \infty $. 
Then there exist $\delta_0 > 0$ and
a simple path $\gamma : [0, \infty) \to \C$  with $\gamma(t) \to \infty$ as $t \to + \infty$ and  the following properties:  
\begin{equation}
 \label{lrw1}
(i) \quad 
\lim_{z \to \infty, z \in \gamma} \frac{u(z)}{\log |z|} = + \infty; \quad (ii) 
\quad \text{if $\lambda > 0$ then} \quad \int_\gamma \exp( - \lambda u(z)) \, |dz| < \infty;
\end{equation}
(iii)
if $z = \gamma (t) $ then  $u(\gamma(s)) \geq \delta _0 u(z) $ for all $s \geq t$. 
\end{thm}
Conclusion (iii) and the fact that $\gamma$ may be chosen to be simple are not stated in \cite[Theorem~1]{LRW},
but both are implicit in the proof. Here $\gamma = \gamma_1 \cup \gamma_2 \cup \ldots $ is constructed in \cite[Section 3]{LRW}
so that, for some fixed $\delta_1 \in (0, 1)$, each $\gamma_k :[k-1, k] \to \C$ is 
a simple path from $a_k \in D_k$ to $a_{k+1} \in \partial D_k$, where $D_k$ is the component  of $\{ z \in \C : \, u(z) < (1-\delta_1)^{-1} u(a_k) \}$ containing $a_1$.
By \cite[(3.2) and (3.3)]{LRW}, the $\gamma_k$ are such that 
$0 < \delta_1 u(a_k) \leq u(z) < (1-\delta_1)^{-1} u(a_k) $ on $\lambda_k = \gamma_k \setminus \{ a_{k+1} \} $ and $u(a_{k+1}) \geq (1-\delta_1)^{-1} u(a_k) > u(a_k) $. 
Hence if $z = \gamma (t) \in \lambda_k   $ then $u( \gamma(s)) \geq \delta_1 u(a_k) \geq \delta_1 (1- \delta_1 ) u(\gamma(t)) $ for all
$s \geq t$. If the whole path $\gamma$ is not simple, take the least $k \geq 2$ such that $\Gamma_k = \gamma_1 \cup 
\ldots \cup \gamma_k$ is not simple. Then there exists a maximal 
$t \in [k-1, k]$ such that $u_k =
\gamma_k(t) $ lies in the compact set $\Gamma_{k-1}$, and $t < k$ since 
$\gamma_k(k) = a_{k+1} \in \partial D_k$. Replacing $\Gamma_k$ by the part of 
$\Gamma_{k-1}$ from $a_1$ to $u_k$, followed by the part of 
$\gamma_k$ from $u_k$ to $a_{k+1}$,  does not affect conclusions (i), (ii) and (iii), and  the argument may then be repeated.

\hfill$\Box$
\vspace{.1in}

Theorem \ref{thmlrw} leads to the following result. 
\begin{prop} 
\label{dtsprop}
Let $N \in \N$ and let
$A$ be a transcendental meromorphic function in the plane such that the inverse function of $A$  has a direct transcendental singularity over $0$.
Then there exist a path $\gamma$ tending to infinity in $\C$ 
and linearly independent solutions $U$, $V$ of 
\begin{equation}
 \label{1}
w'' + A(z) w = 0
\end{equation}
on a simply connected domain containing $\gamma$,  such that $U$ and $V$ satisfy, as $z \to \infty$ on $\gamma$, 
\begin{equation}
 \label{geqn1}
U(z) = z + \frac{O(1)}{z^N} , \quad U'(z) = 1 + \frac{O(1)}{z^N} ,   \quad V(z) = 1 + \frac{O(1)}{z^N} , \quad V'(z) = \frac{O(1)}{z^N} . 
\end{equation}
\end{prop}
To prove Proposition \ref{dtsprop}, observe first that, as in the proof of Theorem \ref{thmds},
there exist a small positive $\delta $ and a
non-empty component
$D $ of  $\{ z \in \C : |A(z)| < \delta \}$ 
such that $A(z) \not = 0$ on $D$,
as well as
a non-constant subharmonic function $u$ on $\C$ given by
$$
u(z) = \log \left| \frac{\delta }{ A(z) } \right|  \quad 
(z \in D ), \quad
u(z) = 0  \quad ( z \not \in D).
$$
Then  $u$ satisfies the hypotheses of Theorem \ref{thmlrw}, by  \cite[Theorem 2.1]{BRS}, and so there exists a path $\gamma : [0, \infty) 
\to D$ as in  conclusions (i), (ii) and (iii). 
In particular, (iii) implies that 
\begin{equation}
 \label{maxAest}
\text{if $z = \gamma (t) $ then  $|A(\gamma(s))| \leq \delta^{1- \delta _0 } |A(z)|^{\delta_0}  $ for all $s \geq t$.}
\end{equation}
Choose a simply connected domain $\Omega$ on which $A$ has no poles, such that $\gamma \subseteq \Omega$. 
By  (\ref{lrw1}) it may be assumed that $|A(t)|^{-1/4} \geq |t|^2 \geq 4$ on $\gamma$, and that 
\begin{equation}
 \label{lrw2}
\int_{\gamma} |t|^2 |A(t)|  \, |dt| \leq 
\int_{\gamma}  |t|^2 |A(t)|^{1/2}  \, |dt| \leq \int_{\gamma}  |A(t)|^{1/4}  \, |dt|  < \frac14 .
\end{equation}

\begin{lem} 
\label{lemgron1}
 Let $v$ be a solution of (\ref{1}) on $\Omega$. Then 
$v(z) = O(|z|)$ 
as $z \to \infty$ on $\gamma$.
\end{lem}
\textit{Proof.} This is a standard argument along the lines of Gronwall's lemma. Let $y_0$ be the starting point of $\gamma$. 
Differentiating twice shows that there exist constants $a_1, b_1$ such that, on $\Omega$, 
$$
v(z)  = a_1 z + b_1  - \int_{y_0}^z (z-t) A(t) v(t) \, dt .
$$
If $\phi(z) = v(z)/z$ is unbounded on $\gamma$ there exist $\zeta_n \to \infty$ on $\gamma$ such that $ \phi( \zeta_n ) \to \infty$ and
$|\phi(t) | \leq | \phi( \zeta_n )|$ 
on the part of $\gamma$ joining $y_0$ to $\zeta_n$. If $n$ is large then  (\ref{lrw2}) delivers a contradiction  via 
$$
| \phi( \zeta_n )| \leq  |a_1| + |b_1| + | \phi( \zeta_n )| \int_{y_0}^z (1+|t|) | t A(t)|  \, |dt| \leq  |a_1| + |b_1| + \frac{| \phi( \zeta_n )|}2 .
$$
\hfill$\Box$
\vspace{.1in}

\begin{lem}
 \label{lemnewLIsolns}
(a) Let $ N \in \N$. Then on $\gamma$ every solution $v_j$ of (\ref{1}) has 
\begin{equation}
 \label{newvrep} 
v_j(z) = \alpha_j z + \beta_j + \int_z^\infty (z-t) A(t) v_j(t) \, dt , \quad \alpha_j, \beta_j \in \C ,
\end{equation}
the integration being from $z$ to infinity along $\gamma$. Moreover,  $v_j$ satisfies, as $z \to \infty$ on $\gamma$, 
\begin{equation}
 \label{bettervest}
v_j(z) - \alpha_j z -  \beta_j = \frac{O(1)}{z^N} , \quad v_j'(z) - \alpha_j  = \frac{O(1)}{z^N} .
\end{equation}
(b) If $v_1, v_2$ are linearly independent solutions  of (\ref{1}) on $\Omega$ then $|\alpha_1| + |\alpha_2| > 0$ in (\ref{newvrep}),
and if $\alpha_2 = 0$ then $\beta_2 \neq 0$. 
\end{lem}
\textit{Proof.} First, (\ref{newvrep}) follows from (\ref{lrw2}) and Lemma \ref{lemgron1}. Next, 
(\ref{lrw1}), (\ref{maxAest}), (\ref{lrw2}), (\ref{newvrep}) and Lemma~\ref{lemgron1} imply that, as $z \to \infty$ on $\gamma$,   
\begin{eqnarray*}
 | v_j(z)- \alpha_j z - \beta_j | &\leq&   |z| \int_{z}^\infty (1+|t|) |A(t)| O( |t| ) \, |dt|\\
 &\leq&   |z| \,  \delta^{(1- \delta _0 )/2} |A(z)|^{\delta_0/2 }  \int_{z}^\infty (1+|t|) |A(t)|^{1/2} O( |t| ) \, |dt| = \frac{O(1)}{z^N} ,\\
|v_j'(z) - \alpha_j| &=&  \left|  \int_{z}^\infty  A(t) v_j(t) \, dt \right| \\
&\leq& \delta^{(1- \delta _0 )/2} |A(z)|^{\delta_0/2 }  \int_{z}^\infty  |A(t)|^{1/2} O( |t| ) \, |dt| = \frac{O(1)}{z^N} .
\end{eqnarray*}
Finally, 
suppose that 
$v_1, v_2$ are linearly independent solutions  of (\ref{1}) on $\Omega$ but the conclusion of 
(b) fails. Then $ v_1(z) v_2'(z) - v_1'(z) v_2(z) \to 0$   as $z \to \infty$ on $\gamma$,  by (\ref{bettervest}), 
contradicting the fact that $W(v_1, v_2)$ is a non-zero constant by Abel's identity. 
\hfill$\Box$
\vspace{.1in}

Now fix linearly independent solutions $v_1, v_2$  of (\ref{1}) on $\Omega$. Then $\alpha_1, \alpha_2$ cannot both vanish in (\ref{newvrep}).
On the other hand, it is possible to ensure that one of $\alpha_1, \alpha_2$ is $0$, by otherwise considering 
$\alpha_2 v_1 - \alpha_1 v_2$. Hence it may be assumed that $\alpha_1 = 1$, while $\alpha_2 =0$ and $\beta_2 = 1$. Now write $U = v_1$ and $V = v_2$, 
so that Lemma \ref{lemnewLIsolns} gives (\ref{geqn1}).
\hfill$\Box$
\vspace{.1in}

\section{Proof of Theorem \ref{thm1aa}}\label{WV}

Assume that $f$ and $S_f$ are as in the hypotheses but that the inverse function of $S_f$  has a direct transcendental singularity over $0$.
Then evidently so has that of $A = S_f/2$, and it is well known 
that (\ref{Schdef}) implies that $f$ is locally the quotient of linearly independent
solutions of (\ref{1}). 
Now Proposition \ref{dtsprop} gives linearly independent solutions $U$, $V$ of (\ref{1}) satisfying (\ref{geqn1}) 
on a path $\gamma $ tending to infinity. 
Moreover, $h = U/V $ has the form $h = T \circ f$, for some  M\"obius transformation $T$, and so $h \in \mathcal{S}$,
whereas $h(z) \sim z$ and $z h'(z)/h(z) = O(1)$ on $\gamma$,  contradicting (\ref{Sineq}).
\hfill$\Box$

\noindent
School of Mathematical Sciences, University of Nottingham, NG7 2RD.\\
james.langley@nottingham.ac.uk


\vspace{.1in}


{\footnotesize
\begin{thebibliography}{1}
\bibitem{Ber4}
W. Bergweiler, Iteration of meromorphic functions, \textit{Bull. Amer. Math.}
Soc. 29 (1993), 151-188.
\bibitem{BE}W. Bergweiler and A. Eremenko, On the singularities of the inverse to a meromorphic function of finite order, 
\textit{Rev. Mat. Iberoamericana} 11 (1995), 355-373.
\bibitem{BEdir}
W. Bergweiler and A. Eremenko, Direct singularities and completely invariant domains of entire functions, 
\textit{Illinois J. Math.} 52 (2008), 243-259. 
\bibitem{BRS}
W. Bergweiler, P.J. Rippon and G.M. Stallard,
Dynamics of meromorphic functions with direct or logarithmic
singularities, \textit{Proc. London Math. Soc.} 97 (2008), 368-400.
\bibitem{CER}
J. Clunie, A. Eremenko and J. Rossi, On equilibrium points of logarithmic
and Newtonian potentials, \textit{J. London Math. Soc.} (2) 47 (1993), 309-320.
\bibitem{elf}
G. Elfving, \"Uber eine Klasse von Riemannschen Fl\"achen und ihre Uniformisierung, \textit{Acta Soc. Sci. Fenn.} 2 (1934) 1-60.
\bibitem{Hay1}
W.K. Hayman, Picard values of meromorphic functions and their derivatives,
\textit{Ann. of Math.} 70 (1959), 9-42.
\bibitem{Hay2} W.K. Hayman, \textit{Meromorphic functions},
Clarendon Press, Oxford, 1964. 
\bibitem{Hay7} W.K. Hayman, \textit{Subharmonic functions Vol. 2}, Academic Press, London, 1989.
\bibitem{Hay9}
W.K. Hayman, \textit{Multivalent functions}, 2nd edition, Cambridge Tracts in
Mathematics 110, Cambridge University Press, Cambridge 1994.
\bibitem{Hin}
J.D. Hinchliffe, The Bergweiler-Eremenko theorem for finite lower order,
\textit{Result. Math.} 43 (2003), 121-128.
\bibitem{La11}
J.K. Langley, The zeros of the first two derivatives of a
meromorphic function, \textit{Proc. Amer. Math. Soc.} 124, no. 8 (1996), 2439-2441.
\bibitem{lajda}
J.K. Langley, Non-real zeros of higher derivatives of real entire functions
of infinite order,
\textit{J. Analyse Math.} 97 (2005), 357-396.
\bibitem{Laschwarzian}
J.K. Langley, The  Schwarzian derivative and the Wiman-Valiron property, \textit{J. d'Analyse Math.} 130 (2016), 71-89. 
\bibitem{LaRo2}
J.K. Langley and John Rossi, Critical points of certain discrete potentials,
\textit{Complex Variables} 49 (2004), 621-637.
\bibitem{LRW} J. Lewis, J. Rossi and A. Weitsman, On the growth of subharmonic functions along paths, \textit{Ark. Mat.} 22 (1983), 104-114.
\bibitem{Nev2}
R. Nevanlinna, \"Uber Riemannsche Fl\"achen mit endlich vielen Windungspunkten,
\textit{ Acta Math.} 58 (1932), 295-373.
\bibitem{Nev}
R. Nevanlinna, \textit{Eindeutige analytische Funktionen,
2. Auflage}, Springer, Berlin, 1953.
\bibitem{sixsmith}D.J. Sixsmith, A new characterisation of the Eremenko-Lyubich class, \textit{J. Analyse Math.} 123 (2014), 95-105. 
\bibitem{Tsuji1}
M. Tsuji, \textit{Potential theory in modern function theory},
Maruzen, Tokyo, 1959.
\end{thebibliography}
}
\end{document}